\documentclass{article}
\usepackage{amsmath, amsthm} 
\usepackage{amssymb} 
\usepackage{bm} 
\begin{document}
%
\long\def \delete #1{} 
\def \degrees {^\circ} 
\def \dist {{\rm dist}}
\def \angle {{\rm angle}}
\def \cont {{\rm cont}}
\def \perp {{\rm perp}}
\def \half {{\scriptstyle{1\over2}}}
\def \inv #1{{#1}^{\rm -1}}
\def \tran #1{{#1}^{\rm T}}
\def \conj #1{\overline{#1}} 
\def \mag #1{\bigl\Vert#1\bigr\Vert}
\def \defn {\equiv} 
\def \floor #1{\lfloor#1\rfloor}
\def \ceil #1{\lceil#1\rceil}
\def \GCD {{\rm GCD}}
\def \QED {\ \square\ } 
\def \fac {\,|\,} 
\def \rnd #1{\langle#1\rangle} 
\def \from {\leftarrow} 
\def \mod #1{\!\!\pmod #1} 
\def \R {\mathbb{R}} 
\def \C {\mathbb{C}} 
\def \H {\mathbb{H}} 
\def \O {\mathbb{O}} 
\def \D {\mathbb{D}} 
\def \N {\mathbb{N}} 
\def \Q {\mathbb{Q}} 
\def \Z {\mathbb{Z}} 
\def \L {\mathbb{L}} 
\def \i {{\,\bf i}} \def \j {{\,\bf j}} \def \k {{\,\bf k}} 
\def \io {{\bm \iota}} 
%
\newenvironment{lrindent}[0]
{\begin{list}{}
  {\setlength{\leftmargin}{\parindent}
  \setlength{\rightmargin}{\parindent}}
  \item[]
}
{\end{list}}
\def\keywords{\vspace{.5em}
{\textit{Keywords}:\ \relax%
}}
\def\endkeywords{\par}
\def\subjclass{\vspace{.5em}
{\textit{AMS Mathematics Subject Classification (2010)}:\ \relax%
}}
\def\endsubjclass{\par}
%
\newtheorem{Thm}{Assertion} 
\newtheorem{Conj}[Thm]{Conjecture} 
\newtheorem{Hypo}[Thm]{Hypothesis} 
\title{Lattice Embedding of Heronian Simplices}
%
%
%
\author{
{W.~Fred Lunnon},
{NUI Maynooth, Rep. of Ireland}
}
\maketitle              
\begin{abstract}
A {\sl rational} triangle has rational edge-lengths and area;
a {\sl rational} tetrahedron has rational faces and volume;
either is {\sl Heronian} when its edge-lengths are integer,
and {\sl proper} when its content is nonzero.
\par
A variant proof is given, via complex number GCD, of the previously
known result that any Heronian triangle
may be embedded in the Cartesian lattice $\Z^2$; it is then shown that,
for a proper triangle, such an embedding is unique modulo lattice isometry;
finally the method is extended via quaternion GCD to tetrahedra in $\Z^3$,
where uniqueness no longer obtains, and embeddings also exist which
are unobtainable by this construction.
\par
The requisite complex and quaternionic number theoretic background is summarised
beforehand. Subsequent sections engage with subsidiary implementation issues:
initial rational embedding, canonical reduction, exhaustive search for
embeddings additional to those yielded via GCD; and illustrative numerical
examples are provided.
\par
A counter-example shows that this approach must fail in higher dimensional
space. Finally alternative approaches by other authors are summarised.
\par
\end{abstract}

%

\section{Introduction: Rational and Heronian Simplices}
\par
A {\sl rational} triangle, tetrahedron, or simplex in {$n$-space} $\R^n$ is
a polytope with $n+1$ vertices, and rational edge lengths, face areas,
solid volumes, etc. When the edges have integer length, the
simplex is {\sl Heronian}; when the GCD (greatest common divisor)
of those lengths equals unity, it is {\sl primitive} (Heronian).
When its $n$-dimensional content is nonzero, it is {\sl proper} ---
we generally have only proper cases in mind, despite which much of the
following may easily be seen to apply more widely.
Obviously any rational simplex may be dilated to a Heronian simplex. 
\par
Rational simplices are of interest as test-cases since many of their subsidiary
elements are likewise rational --- for triangles, the
edges, altitudes, medians, centroid, orthocentre, incentre, circumcentre, 
angle sines and cosines, etc.
\par
Our purpose here is to show that for $n = 2,3$ the simplex in $\R^n$
may be posed first {\sl axially} in $\Q^n$ (which is elementary), and finally
{\sl embedded} in the lattice $\Z^n$; that is, congruently so that its
vertices have Cartesian coordinates not just rational, but integer.
For triangles the result is only about 10 years old \cite{[Yiu01]};
we propose a construction based on complex-number GCD, and related to the
approach of \cite{[Fri01]}. 
\par
For tetrahedra the paradigm must be extended to embrace quaternion GCD,
where commutativity and unique factorisation are no longer available:
the correspondence between the two arguments is illuminating,
and we emphasise the analogies and contrasts by presenting the main theorems
in sections \ref{[embedR2]}, \ref{[embedR3]} in copycat format.
The level of detail in their proofs might be felt excessive, to the point
of obscuring the wood for the trees: it was motivated by a desire to avoid
elementary blunders which have compromised numerous earlier versions.
\par
Incidentally, an obscure but significant step in the argument involves first
establishing that the preliminary rational pose (on which GCD is to act)
has denominators with prime factors $p\equiv 1\mod 4$ only.
\par
It is natural to enquire how these results might generalise to higher
dimensions. There are currently several obstructions: the quaternion
representation employed cannot be extended further; even if it
could, the theorem in the form given here would no longer hold; and in any
case, no examples of Heronian pentatopes are known to which any such theorem
might be applicable.
\par
The contribution to this project from Warren D.~Smith, in the form of copious
literature references, enthusiastic criticism, and awesome energy,
is gratefully acknowledged. Some motivating ideas, particularly behind the
proof of the embedding theorem, originated with Michael Reid. Finally
Susan Marshall and Alexander Perlis were instrumental in uncovering yet
another embarrassing error in an earlier version of this proof.
\par
Despite having spent much time contemplating the matter, this author still
finds the connection revealed, between the number theoretic concept of GCD
and the geometric concept of lattice embedding, both surprising and mysterious.
\par

\section{Edges, Triangles, Tetrahedra and their Content} \label{[arevol]}
\par
The Hero formula giving the area $d$ of a triangle with edge lengths $u,v,w$ as
\begin{align} \label{[Hero]}
(4 d)^2 = (u+v+w)(u+v-w)(u-v+w)(-u+v+w)
\end{align}
is case $n = 2$ of the {Cayley-Menger} formula for the content of a
simplex in $\R^n$, which deserves to be better known \cite{[Som30]}:
\begin{align} 
(4 d)^2 = 
\left|\begin{array}{c c c c} \label{[caymen2]}
  0 & 1 & 1 & 1 \cr
  1 & 0 & u^2 & v^2 \cr
  1 & u^2 & 0 & w^2 \cr
  1 & v^2 & w^2 & 0 
\end{array}\right| .
\end{align}
\par
For $n = 3$ the volume $e$ of the tetrahedron with edge lengths $u,\ldots,z$
is given by
\begin{align} 
2(12 e)^2 = 
\left|\begin{array}{c c c c c} \label{[caymen3]} 
  0 & 1 & 1 & 1 & 1 \cr
  1 & 0 & u^2 & v^2 & x^2 \cr
  1 & u^2 & 0 & w^2 & y^2 \cr
  1 & v^2 & w^2 & 0 & z^2 \cr
  1 & x^2 & y^2 & z^2 & 0 
\end{array}\right| ;
\end{align}
Since these coefficients turn out to be even, $(12 e)^2$ is again a homogeneous
symmetric polynomial with integer coefficients, but now cubic in the edge
squares.
\par
Even less familiar are rational expressions for angular quantities:
for a triangle,
\[ \cot p/2 = (s-w)s/d \]
etc., where $p$ denotes the angle at vertex $P$. The rationality of these implies that all trigonometric functions of the angles are also rational,
a property christened {\sl geodetic} by Conway. Relations between angular
and linear quantities are pursued in section \ref{[1mod4]}, where they find
an unexpectedly arcane application.
\par
Heronian area and volume have interesting divisibility properties, which we
regretfully must refrain from discussing.
\par
To specify a particular free tetrahedron, a correspondence must be established
between vertices and edges: denoting vertices $P,Q,R,S$, edge-lengths will be
specified by a hexad of integers presented in sequential order 
\[ [QP,RP,RQ,SP,SQ,SR]. \]
\par
The amount of data involved is reduced by restriction to {\sl primitive}
cases, in which the GCD (greatest common divisor) of the edges equals unity;
also by reduction to {\sl canonical form}, with vertices permuted so that
the vector of edge-lengths represents the largest possible `number', where
individual components are regarded as `digits'. [A humbling exercise
in algorithmic design is to sort such a hexad into canonical form efficiently,
avoiding inspection of all 24 tetrahedral symmetries.]
\par
Where rational Cartesian coordinates are involved, we frequently refer to
the LCD (least common denominator) of their components: that is, the LCM
(least common multiple) of the denominators after each component has been
reduced so that GCD of numerator and denominator equal unity.
When the coordinates are instead projective, this concept is reformulated as
the {\sl primitive} reduction of a rational vector with scalar component unity:
multiplying by the LCM of the denominators, then dividing by the GCD of the
resulting integers, causes the LCD to appear as scalar component.
\par

\section{Representation of $\R^2$ and $\Z^2$ via Complex Numbers}
\par
With $x,y$ real in $\R$, denote a complex number in $\C = \R[\io]$ by
\[ Z = x + \io y ,\]
its {\sl norm} by
\[ \mag{Z} \equiv x^2 + y^2 ,\]
its {\sl conjugate} by
\[ \conj{Z} \equiv x - \io y .\]
\par
Any $P = x + \io y$ represents a finite point with Cartesian coordinate
$(x,y)$. The representation is {\sl inhomogeneous}: scalar multiples of $P$
represent distinct points.
The origin point is represented by $X = 0$.
The distance squared between points $P, Q$ equals $\mag{P - Q}$.
\par
Any {\sl unit} $U = \cos t + \io \sin t$ represents rotation through
angle $t$ around the origin; non-unit $U$ with $\mag{U} \ne 1$
do not represent rotations. The identity rotation is represented by $U = 1$.
Composition of rotations is represented by commutative complex product.
\par
The transformation
\[ P \to U P \quad{\rm where}\quad U = \cos t + \io \sin t \]
rotates the point $P = x + \io y$ through angle $t$.
\par
A {\sl (Gaussian) integer} lies in $\Z[\io]$, its components rational integers.
A complex {\sl prime} is either a rational prime $p \equiv -1\mod 4$, or of
form $x \pm \io y$ where $x^2+y^2 = p$ with $p = 2$ or $p \equiv 1\mod 4$.
Factorisation into primes is essentially unique modulo unit {\sl association}:
multiplication of individual primes by $\pm 1, \pm \io$.
\par
$\Z[\io]$ is a Euclidean domain: given $X, Y \in \Z[\io]$,
their greatest common divisor $\GCD(X, Y)$ is unique modulo association.
Algorithmic details of GCD are discussed in section \ref{[GCDalg]}.
For further detail see \cite{[Con03]} chap.~2.
\par
Unit integers represent rotations through $k\pi/2$ around the origin.
The full group of lattice isometries is generated by these 4 rotations,
reflections in axes (complex conjugacy etc), and translations (integer addition).

\section{Lattice Embedding Theorem in $\R^2$} \label{[embedR2]}
\par
\begin{Thm} \label{[poseone2]}
Suppose all but one of a set $\{Q\}$ of $k+1$ points are already embedded in
$\Z^2$; the other point $R$ is rational in $\Q^2$ with LCD $r$;
and the squares of all their point-to-point distances lie in $\Z$.
Then rotation via complex product $Q \to \conj{X}^2 Q/r$ where $X =\GCD(r R, r)$
embeds the entire set congruently in $\Z^2$.
\end{Thm}
Given the first assertion, translating some point of the set to the origin
and applying induction on the number $k$ of points, immediately
\begin{Thm} \label{[poseall2]}
If a set of points lies in $\Q^2$, and their squared distances lie in $\Z$,
then the entire set may be embedded congruently in $\Z^2$.
\end{Thm}
And taking the set to be the vertices of a Heronian triangle in some
rational pose (section \ref{[axialpose]}), with edge-lengths themselves integer,
\begin{Thm} \label{[embed2]}
Any Heronian triangle may be embedded congruently in $\Z^2$.
\end{Thm}
\par
\begin{proof} (Assertion \ref{[poseone2]}) 
Let the points $P,Q,\ldots$ lie in $\Z^2$, except for $R$ in $\Q^2$
with Cartesian coordinate $(x/r, y/r)$ and component LCD $r$.
Represented in $\C$, all lie in $\Z[\io]$ except for $R$;
now dilating the entire set by $r$ yields also $R'' = r R \in \Z[\io]$.
\par
Some lattice point $P$ of the set may be translated to the origin $P = 0$;
now if $k = 0$ the result becomes trivial. Otherwise since all distances squared
$\mag{P - Q}$ are integer, all dilated points have magnitude divisible by $r^2$.
\par
Set $X = \GCD(R'', r)$, for which evidently $\mag{X}\fac r$.
Every rational prime factor $p$ of $r$ corresponds to some properly complex
prime factor $W$ of $X$; otherwise we should have $p \fac R''$,
which was primitive. Therefore in fact
\[\mag{X} = r = X \conj{X} = \conj{X} X .\]
\par
Firstly, the rational point $R$ is transformed by $X$ to the lattice.
For $r^2 \fac R''$, and each factor $W$ must occur twice over in $R''$ --- 
otherwise both $W$ and $\conj{W}$ would occur, and again $p \fac R''$.
Therefore $X^2 | R''$; and $R$ is transformed to 
\[ \conj{X}^2 R''/r^2  =  \conj{X}^2 X^2 R'/r^2  =  R' \]
for some $R' = \conj{X}^2 R/r \in \Z[\io]$.
\par
Secondly, each lattice point $Q$ is transformed by $X$ to the lattice.
For
\[ r^2\ \fac\ \mag{r Q - R''} = (r Q - R'') \conj{(r Q - R'')} \]
which being scalar is transformed to (and equals)
\[ \mag{\conj{X}^2 r Q - \conj{X}^2 R''}/r^2 = 
(\conj{X}^2 r Q - \conj{X}^2 R'')(X^2 r \conj{Q} - X^2 \conj{R''})/r^2; \]
hence
\[ r^4\ \fac\ (\conj{X}^2 r Q - \conj{X}^2 R'')(X^2 r \conj{Q} - X^2 \conj{R''}). \]
Now for $r \in \Z$, $Y \in \Z[\io]$ we have $r\fac Y$ iff $r\fac \conj{Y}$;
so via unique factorisation
\[ r^2\ \fac\ \conj{X}^2 r Q - \conj{X}^2 R''; \quad{\rm also}\quad 
r^2\ \fac\ \conj{X}^2 R'' ,\]
and adding, finally
\[ r^2\ \fac\ r \conj{X}^2 Q .\]
Therefore $Q$ is transformed to $Q'= \conj{X}^2 Q/r \in \Z[\io]$.
\end{proof} 
\par
The construction is in fact reversible, although this seems difficult to
establish directly.
The following argument relies heavily on unique factorisation.
\begin{Thm} \label{[unique2]}
Lattice embedding of a proper Heronian triangle is unique modulo lattice
isometry (unit association, integer addition, conjugacy).
\end{Thm}
\begin{proof}
Suppose $P,Q,R$ and $P,S,T$ are lattice embeddings into $\Z[\io]$
from vertices of the same free triangle, translated and reflected to 
share the origin $P = 0$ and to have anticlockwise sense.
\par
Let $Q' = Q/\GCD(Q, S)$, $S' = S/\GCD(Q, S)$, and set $U = S/Q = S'/Q'$;
the isometry $U$ rotates $P,Q,R$ to $P,S,T$.
Since $\mag{Q} = \mag{S}$ we have $\mag{Q'} = \mag{S'} = q$ say.
\par
If $q = 1$ then $U \in \Z[\io]$ is a lattice isometry, 
so the embeddings are essentially equivalent.
Otherwise, each rational prime factor of $q$ splits into conjugate primes
in $Q'$ and $S'$, which were coprime;
so $Q' = \conj{S'}$ and $U = S'/Q' = {S'}^2/q$, and similarly $U = {T'}^2/q$.
But now $S' = \pm T'$, so the triangle was improper.
\end{proof}
\par

\section{Representation of $\R^3$ and $\Z^3$ via Quaternions}
\par
With $s,p,q,r$ real in $\R$, 
\[ X = s + p \i + q \j + r \k \]
denotes a quaternion in $\H = \R[\i,\j,\k]$;
\[ \mag{X} \equiv s^2 + p^2 + q^2 + r^2 \]
denotes its {\sl norm};
\[ \conj{X} \equiv s - p \i - q \j - r \k \]
denotes its {\sl conjugate}.
\par
In general $X \ne 0$ represents the finite point at distances $p/s, q/s, r/s$
from the Cartesian coordinate planes. The representation is {\sl homogeneous}:
any nonzero scalar multiple of $X$ represents the same point.
The origin point is represented by $X = 1$ (or any nonzero scalar).
A {\sl pure} quaternion with scalar component $s = 0$ represents a point at
projective infinity.
The distance squared between points $P, Q$ equals $\mag{P - Q}$
only when both have been {\sl co-normalised}, with scalar parts unity.
\par
Dually $X \ne 0$ represents rotation through angle $t$ around a line
meeting the origin, where $s = \cot t/2$ and the direction cosines
(of angles made by the line with the coordinate planes) are $p, q, r$.
The representation is homogeneous: any nonzero scalar multiple of $X$
represents the same rotation.
The identity rotation is represented by $X = 1$ (or any nonzero scalar).
Composition of rotations is represented by quaternion product,
which is in general non-commutative, unless its arguments differ
only in their scalar component.
\par
In vector terms, given 3-vector $u$, unit 3-vector $v$, and angle $t$,
the {\sl transformation} (group-theoretic conjugation)
\[ P \to \conj{X} P X \quad{\rm where}\quad X = \cot t/2 + v \]
rotates the point $P = 1 + u$ through angle $t$ about axis $v$.
\par
Interpreted as rotations, quaternions are analogous to conventional
{\sl polar} or {\sl contravariant} vectors; as points, to {\sl axial}
or {\sl covariant} vectors: and this duality accounts for their varying
normalisations. Transformation by rotation is the only situation
in which it is legitimate to multiply a point.
\par
Symmetries of the Cartesian lattice are generated by translations
(integer additions), reflection in the origin (conjugacy),
and 24 rotations fixing the unit cube via these transformations:
\begin{align*}
  &1 &\ ({\rm identity\ }) \times 1; \cr
  &\i,\j,\k &\ (\pi {\rm\ around\ axis\ })\times 3; \cr
  &(1\pm\i)/\sqrt2,\ldots &\ (\pi/2  {\rm\ around\ axis\ })\times 6; \cr
  &(\i\pm\j)/\sqrt2,\ldots &\ (\pi/2 {\rm\ around\ edge\ perpendicular\ })\times 6; \cr
  &(1\pm\i\pm\j\pm\k)/2 &\ (2\pi/3 {\rm\ around\ diagonal\ })\times 8. 
\end{align*}
\par

A {\sl (Lipschitzian) integer} $X \in \L = \Z[\i,\j,\k]$ has
rational integer components. It is {\sl primitive} when the GCD of its
components equals unity. Rational points are representable by $X \in \L$, 
lattice points by $X \in \L$ with unit scalar, and conversely.
\par
Quaternion number theory is complicated by non-commutativity of multiplication and non-unique factorisation.
A {\sl prime} $X \in \L$ has prime norm $\mag{X} = p$,
while $X \notin \Z$ itself; a rational prime may factorise into quaternion
primes in several essentially distinct ways.
Factorisation is discussed in detail in \cite{[Con03]} chap.~5.
\par
Integer quaternions are almost a Euclidean domain: provided at least one of $Y,
Z \in \L$ has odd norm, their left- and right-hand greatest common divisors
\[ \GCD_L(Y, Z),\quad \GCD_R(Y, Z) \]
are guaranteed to exist, though in general distinct, and are unique modulo 
{\sl association}: product with units $\pm 1, \pm \i, \pm \j, \pm\k$
on the right, left respectively. $X = \GCD_L(Y, Z)$ has the property that
$X$ divides $Y,Z$ on the left --- $Y = X U$ and $Z = X V$ for some $U,V \in \L$
--- and for any $W$ dividing $Y,Z$ on the left, $W$ divides $X$ on the left;
analogously for $\GCD_R$.
\par
Factorisation of an arbitrary $X \in \L$ is no longer unique.
In the first place, within some factorisation of $X$ the product $Y Z$ of
any two factors may trivially be replaced by $(Y J)(\conj{J} Z)$,
with $J \in \L$ a unit: this process is denoted {\sl unit migration},
and as in real and complex factorisation we consider identical factorisations
equivalent modulo unit migration.
\par
Now suppose the norm of $X$ factorised into rational primes, some possibly
equal, as
\[ \mag{X} = p_1\ldots p_i\ldots \]
in some arbitrarily specified order. Then provided $X$ is primitive with odd
norm, $Y = \GCD_L(X, p_1)$ is the unique left factor of $X$ with norm $p$
modulo association; now by induction on $i$, the factorisation of $X$ proceeds
uniquely modulo unit migration. In summary:
\begin{Thm} \label{[unifactor]}
Given any specified order for its factor norms, an odd primitive
$X \in \L$ factorises uniquely into primes modulo unit migration.
\end{Thm}
\par
The existence of the GCD may
be inferred from the termination of the algorithm in section \ref{[GCDalg]};
see also \cite{[Dic22]} theorem~3.
Provided only $q \fac Q$ and $q \in \Z$, then via \cite{[Dic24]} theorem~2
there do still exist $X, W \in \L$ such that $q = X \conj{X}$ and $Q = X W$;
but only when $q$ is odd can $X = \GCD_L(Q, q)$ be guaranteed unique.
Analogous observations apply on the right-hand side.
\par

\section{Algorithms for GCD} \label{[GCDalg]}
\par
Define real and quaternion {\sl rounding} to integer via
\[ \rnd{x} \equiv \floor{x + 1/2} \in \Z \quad{\rm for\ } x \in \R ,\]
\[ \rnd{X} = \rnd{s + p \i + q \j + r \k} \equiv \rnd{s} + \rnd{p} \i + \rnd{q} \j + \rnd{r} \k \in \Z[\i, \j, \k] \quad{\rm for\ } X \in \H ;\]
then left and right {\sl remainder} functions are defined by
\begin{align*}
  X {\rm\,mod_L\,} Y &\equiv X - Y \rnd{\conj{Y} X/\mag{Y}}, \cr
  X {\rm\,mod_R\,} Y &\equiv X - \rnd{X \conj{Y}/\mag{Y}} Y; 
\end{align*}
where $\conj{Y}/\mag{Y} = \inv{Y}$ equals the commutative multiplicative
inverse of $Y$.
\par
Now the Euclidean algorithm can be deployed, enhanced to detect cases where no
GCD exists:
\delete{
\begin{align*}
&\quad \GCD_L(Y, Z) \from \\
&\quad\quad {\bf while\ } Z \ne 0 {\bf\ do} \\
&\quad\quad\quad X \from Y;\ Y \from Z;\ Z \from X{\rm\,mod_L\,}Y; \\
&\quad\quad\quad {\bf if\ } \mag{Z} \ge \mag{Y} {\bf\ then\ } {\rm ABORT} {\bf\ f{}i\ od}; \\
&\quad\quad {\bf return\ }Y; 
\end{align*}
}
\begin{equation*}
\boxed{\begin{array}{l}
\GCD_L(Y, Z) \from \\
\quad {\bf while\ } Z \ne 0 {\bf\ do} \\
\quad\quad X \from Y;\ Y \from Z;\ Z \from X{\rm\,mod_L\,}Y; \\
\quad\quad {\bf if\ } \mag{Z} \ge \mag{Y} {\bf\ then\ } {\rm ABORT} {\bf\ f{}i\ od}; \\
\quad {\bf return\ }Y; 
\end{array}}
\end{equation*}
For $\GCD_R(Y, Z)$, substitute ${\rm\,mod_R\,}$ for ${\rm\,mod_L\,}$ above. 
\par
This procedure is guaranteed to terminate provided at least one of $\mag{Y}$
or $\mag{Z}$ is odd, since the only circumstance under which the nonnegative
integer $\mag{Z}$ can fail to further reduce $\mag{Y}$ is when all four
fractional parts equal exactly one half.
But for example when $Y = 2$, $Z = 1 + \i + \j + \k$ (both with norm 4),
the distinct primes $1 + \i,\ 1 + \j,\ 1 + \k$ all divide both,
while no product of all three (with norm 8) can possibly divide either;
so no common divisor $X$ can be divisible by every divisor of both.
\par
GCD in $\C$ follows lines similar to GCD in $\H$, though simplified by
commutativity, just two components, and existence for all arguments.
\par

\section{Rational Embedding of Heronian Tetrahedra} \label{[axialpose]}
\begin{Thm} \label{[ratpose]}
In order that a simplex in $\R^n$ be rationally embeddable in $\Q^n$,
it suffices that the edge-lengths have rational squares, and there exists
some flag (incrementally ascending sequence of vertex subsets) with only
rational contents.
\end{Thm}
The result extends in a natural fashion to polytopes with more than $n+1$
vertices, see \cite{[Fri01]} Prop.~3.3; here we are principally concerned
only with the case $n = 3$.
\par
\begin{proof} 
Given an arbitrary free tetrahedron with edge lengths $u,v,w,x,y,z$ to be
posed in $\Q^3$, its vertices may be posed in axial subspaces of increasing
dimension thus:
\[ P = 1,\quad Q = 1 + Q_1 \i,\quad R = 1 + R_1 \i + R_2 \j,
  \quad S = 1 + S_1 \i + S_2 \j + S_3 \k . \]
Elimination of the 6 equations $\mag{P - Q} = u$ etc. displays the vertex coordinates as {\sl rational} functions of edge length $u$ of $PQ$, face area
$d$ of $PQR$, volume $e$ of $PQRS$, and squares of other edge lengths:
\begin{align*}
  Q_1 &= u, \cr
  R_1 &= (v^2 - w^2 + Q_1{}^2)/2\,Q_1, \cr
  R_2 &= 2\,d/u, \cr
  S_1 &= (x^2 - y^2 + Q_1{}^2)/2\,Q_1, \cr
  S_2 &= (x^2 - z^2 + R_1{}^2 + R_2{}^2 - 2\,S_1\,R_1)/2\,R_2, \cr
  S_3 &= 3\,e/d . 
\end{align*}
\end{proof}
\par
Each vertex may then be represented by some primitive in $\L$,
with scalar equal to the LCD of its coordinate components,
for submission to the lattice embedding procedure in section \ref{[embedR3]}.
[The expressions are only valid for $u,d$ nonzero:
a more elaborate algorithm is required for general improper cases.]
\par
The denominators of these rationals satisfy a constraint which will later find
unexpected application in theorem \ref{[poseone3]}.
\begin{Thm} \label{[1mod4]}
The denominators of coordinate components of vertices of a Heronian tetrahedron
in axial pose are products of primes $p\equiv 1\mod 4$.
\end{Thm}
\par
\begin{proof}
Borrowing a technique from \cite{[Che29]}, the coordinate components
are expressed as homogeneous rationals in terms of half-angle cotangents:
at triangular vertices in $\R^2$, and at dihedral edges in $\R^3$. These
cotangents are rationals $f'/f'' \in \Q$ with $\GCD(f', f'') = 1$;
when substituted into the coordinates, they give rise to denominators
with factors of form $f'^2 + f''^2$.
\par
A theorem from elementary number theory due to Euler \cite{[Fer10]}
guarantees that $f'^2 + f''^2$ factorises into primes $p = 4k + 1$, possibly
with a single extra factor 2 when both $f', f''$ are odd. It then remains
only to check that these surplus 2's cancel with the numerator in each
component.
\par
So consider the $\R^3$ vertex $S$ above, which subsumes all other cases.
Denote by $k = 2\,c/u$ the height of face $PSQ$; at vertex $P$ we find $g$
in terms of $S$ etc.
\begin{align*}
  \cos P &= S_1/x, \quad \sin P = k/x, \cr
  g = \cot P/2 &= (1 + \cos P)/\sin P = (x + S_1)/k; 
\end{align*}
and similarly $h$ at dihedral $U$ along edge $PQ$
\begin{align*}
  \cos U &= S_2/k, \quad \sin U = S_3/k, \cr
  h = \cot U/2 &= (1 + \cos U)/\sin U = (2\,c/u + S_2)/S_3; 
\end{align*}
then inverting to get $S$ in terms of $g,h$ etc.
\begin{align*}
S = 1 &+ x \cos P\i + x \sin P \cos U \j + x \sin P \sin U \k \cr
  = 1 &+ (g^2-1)/(g^2+1)\cdot x \i \cr
      &+ 2\,g(h^2-1)/(g^2+1)(h^2+1)\cdot x \j \cr
      &+ 4\,g h/(g^2+1)(h^2+1)\cdot x \k \cr
  = 1 &+ (g'^2 - g''^2)/(g'^2 + g''^2)\cdot x \i \cr
      &+ 2\,g'g''(h'^2 - h''^2)/(g'^2 + g''^2)(h'^2 + h''^2)\cdot x \j \cr
      &+ 4\,g' g'' h' h''/(g'^2 + g''^2)(h'^2 + h''^2)\cdot x \k ,
\end{align*}
where $g = g'/g''$ and $h = h'/h''$ are positive reduced fractions.
If it is even, $(f'^2 + f''^2)\equiv 2 \mod 4$ and $(f'^2 - f''^2)$
is also even; so all 2's cancel from the denominators as required.
\end{proof}
\par

\section{Lattice Embedding Theorem in $\R^3$} \label{[embedR3]}
\par
\begin{Thm} \label{[poseone3]}
Suppose all but one of a set $\{Q\}$ of $k+1$ points are already embedded in
$\Z^3$; the other point $S$ is rational in $\Q^3$ with LCD $s$;
and the squares of all their point-to-point distances lie in $\Z$.
Then rotation via quaternion transformation $Q \to \conj{X} Q X$ where
$X =\GCD_L(S, s)$ embeds the entire set congruently in $\Z^3$.
\end{Thm}
Given the first assertion, translating some point of the set to the origin,
and applying induction on the number $k$ of points, immediately
\begin{Thm} \label{[poseall3]}
If a set of points lies in $\Q^3$, and their squared distances lie in $\Z$,
then the entire set may be embedded congruently in $\Z^3$.
\end{Thm}
And taking the set to be the vertices of a Heronian tetrahedron in some
rational pose (section \ref{[axialpose]}) with edge-lengths themselves integer,
\begin{Thm} \label{[embed3]}
Any Heronian tetrahedron may be embedded congruently in $\Z^3$.
\end{Thm}
\par
\begin{proof} (Assertion \ref{[poseone3]}) 
Let the points $P,Q,\ldots$ lie in $\Z^3$, except for $S$ in $\Q^3$
with Cartesian coordinate $(x/s, y/s, z/s)$ and component LCD $s$.
Represented in $\L$, every point is primitive;
and each has scalar component unity, excepting $S = s + x\i + y\j + z\k$.
By assertion \ref{[1mod4]}, $s$ may be assumed odd.
\par
Some lattice point $P$ of the set may be translated to the origin $P = 1$;
now if $k = 0$ the result becomes trivial. Otherwise since all distances squared
$\mag{P - Q}$ are integer, $s^2 \fac \mag{S - s P}$, so $s^2 \fac \mag{S}$.
\par
By section \ref{[quatint]} there exists $X = \GCD_L(S, s)$, and evidently
$\mag{X}\fac s$. Every rational prime factor $p$ of $s$ corresponds to some
properly quaternion prime factor of $X$; otherwise we should have $p \fac S$,
which was primitive. Therefore in fact
\[\mag{X} = s = X \conj{X} = \conj{X} X .\]
\par
Firstly, the rational point $S$ is transformed by $X$ to the lattice.
For
\[X = \GCD_L(S, s) = \GCD_L(2s - S, s) = \GCD_L(\conj{S}, 
\conj{s}) = \conj{\GCD_R(S, s)},\]
so
\[\GCD_L(S, s) = X \quad{\rm and}\quad \GCD_R(S, s) = \conj{X} .\]
\par
Choosing some fixed order for the rational prime factors of $s$,
imposing it on $X$, on $S$ at the left, and reversed on $S$ at the right,
yields via Assertion \ref{[unifactor]} the factorisation $S = X S' \conj{X}$,
with $S' \in \L$. [Its left and right factors are necessarily disjoint,
since $\mag{X \conj{X}} = s^2$ and $s^2 \fac \mag{S}$.]
\par
Therefore $S$ is transformed to
\[\conj{X} S X = \conj{X} X S' \conj{X} X = s^2 S';\]
the scalar component of each side equals $s^2$,
so $S'$ has unit scalar, and {\sl qua} point $S' \in \Z^3$.
\par
Secondly, each lattice point $Q$ is transformed by $X$ to the lattice.
For note that $s Q - S = -\conj{(s Q - S)}$ is pure, so 
\[ s^2\ \fac\ \mag{s Q - S} = -(s Q - S)^2 \]
which is transformed to (scaled by $s^2$)
\[ s^4\ \fac\ \mag{\conj{X} s Q X - \conj{X} S X} = -(\conj{X} s Q X - \conj{X} S X)^2; \]
so
\[s^2\ \fac\ s \conj{X} Q X - \conj{X} S X, \quad s^2\ \fac\ \conj{X} S X, 
\quad s^2\ \fac\ s \conj{X} Q X \]
via commutativity of scalars, and adding.
That is, $Q$ is transformed to $\conj{X} Q X = s Q'$, where $Q'\in \L$
with unit scalar as before, so {\sl qua} point $Q' \in \Z^3$.
\end{proof} 
\par
In contrast to dimension~2, permuting vertices may yield embeddings
which are essentially distinct, in the sense defined in section \ref{[canon]}.
Consider the penultimate stage of the embedding procedure,
where $P,Q,R \in \Z^3$ and only $S$ remains to be located: since $P,Q,R$ are
essentially unique via assertion \ref{[unique2]}, at most 4 distinct
configurations are possible, depending upon which vertex was initially
specified as $S$. Therefore
\begin{Thm} \label{[4poses3]}
The number of essentially distinct embeddings of a Heronian tetrahedron into
$\Z^3$, as constructed via the GCD procedure, is at most 4.
\end{Thm}
However, it will emerge in section \ref{[canon]} that there are embeddings
which are not constructible in this fashion. Whether there is a bound on the
total number of distinct embeddings is not known: modulo lattice isometry,
it would have to be at least 36.
\par
Finally, complex numbers turn out to have been something of a waste of space
all along:
\begin{Thm} \label{[emb2emb3]}
Embedding in $\Z^2$ is a special case of embedding in $\Z^3$.
\end{Thm}
\begin{proof}
The claim is equivalent to the following.
Given a rational point $R = 1 + x/r \i + y/r \j$ in $\Q^2$,
the quaternion GCD rotates it to some lattice point remaining in $\Z^2$;
that is, $X = \GCD_L(r + x \i + y \j, r) = f + g \k$ for some $f,g \in \Z$.
To show this, 
\[ \GCD_L(r + x \i + y \j, \,r) = \GCD_L(r, \,(x \bmod r)\i + (y \bmod r)\j) \]
remaindering modulo $r$;
\[ = \GCD_L(r, \,x \bmod r + (y \bmod r)\k) \]
via right product of second argument by $-\i$ (units are coprime);
\[ = \GCD(r, \,x \bmod r + \io\,y \bmod r) \]
setting $\k \equiv \io$, its polar complex equivalent;
\[ = f + \io\,g ,\]
say.
\end{proof}
\par
Notice how (polar) rotations $\cos t + \io \sin t \in \C$
correspond to $\cos t + \sin t \cdot \k \in \H$; whereas (axial)
points $x + \io y \in \C$ correspond to $1 + x \i + y \j \in \H$.
\par

\section{Essentially Distinct Embeddings; Canonical Reduction} \label{[canon]}
\par
Here we consider only tetrahedra in $\R^3$; generalisation to simplices
in $\R^n$ would be routine if required.
\par
While there are infinitely many lattice embeddings of any Heronian tetrahedron,
factoring out translations by fixing one vertex at the origin reduces this
to a finite number. [For technical reasons however, it proves more convenient
instead to further translate so that all coordinates are positive, with
each component zero for at least one vertex.]
\par
Many of these fall into sets of (in general) 48 lattice isomorphs of one another: to eliminate these trivial variations from consideration,
define the {\sl weak canonical} form to be the earliest isomorph, modulo lattice
symmetries, in numerical order --- where coordinate components taken in the
natural ordering are regarded as individual `digits' in a 12-digit `number'.
\par
However even with this refinement embeddings may proliferate,
particularly (and paradoxically) for symmetric cases. They may further be
reduced by a potential factor 4 via permutation of vertices;
so define the {\sl strong canonical} form to be the earliest weak
form modulo vertex permutation.
[The relationship between input edge lengths and vertices is disrupted,
but can straightforwardly be recovered when required.]
\par
The symmetry types occurring in Heronian tetrahedra are of interest in this
connection. The table records number of isomorphs for each group,
along with the corresponding pattern of equal edges, and a (permuted)
example where applicable.
\par
\smallskip
\begin{tabular*}{0.9\textwidth}{@{\extracolsep{\fill}}| c c c c |}
\hline
scalene & 1 & $[u,v,w,x,y,z]$ & [117,84,51,80,53,52] \cr
semi-isosceles & 2 & $[u,v,v,x,x,z]$ & [680,680,208,615,185,185] \cr
isosceles & 4 & $[u,v,v,v,v,z]$ & [1073,1073,990,896,1073,1073] \cr
semi-isohedral & 2 & $[u,v,w,w,v,z]$ & [990,901,793,793,901,308] \cr
isohedral & 4 & $[u,v,w,w,v,u]$ & [203,195,148,148,195,203] \cr
\hline
isohedral-isosceles & 8 & $[u,v,v,v,v,u]$ & NONE \cr
equilateral & 6 & $[u,u,u,x,x,x]$ & NONE \cr
regular & 24 & $[u,u,u,u,u,u]$ & NONE \cr
\hline
\end{tabular*}
\smallskip
\hfill\break
The last two cases above are trivially impossible, since a unit equilateral
triangle has irrational area $\sqrt{3}/2$; the 8-symmetric case is dispatched
at \cite{[Buc92]} sect.~3 case 2(ii).
\par
\delete{
Given the set of all embeddings of a given case for which some chosen vertex
occupies the origin, consider the set of rotations taking one to another.
It is natural to enquire about the structure of this latter set,
in particular, whether it constitutes a group under composition.
But the answer in general is negative, see \cite{[Con03]} chap.~3:
any such group would be finite and contain the octahedral group as a factor,
whence it must actually be the octahedral group, and in that case there is
essentially just one embedding.
\par}
The first Heronian triangle for which the lattice embedding is not also
an axial pose is
\[[30,29,5];\]
the first tetrahedron with some non-axial embedding, also the first without an
axial embedding, is
\[[160,153,25,120,56,39];\]
the first with some embedding not constructible via GCD from any axial pose is
\[ [888,875,533,533,875,888];\]
the first with 2,3,4 distinct embeddings via GCD are respectively
\begin{align*}
[1073,975,448,495,952,840],  \cr
[1360,1092,548,975,865,663], \cr
[45100,43911,6929,34476,40544,36975].
\end{align*}
The particularly prolix case
\[ [8484,6625,6409,6409,6625,8484],\]
yields 9(36) strongly(weakly) distinct embeddings, versus just 1(4) via GCD.
\par

\section{Exhaustive Search for Lattice Embeddings} \label{[search]}
\par
The quaternion GCD-based construction establishes that some lattice
embedding exists for every Heronian tetrahedron, and is
furthermore fast --- far more time is spent in reduction to canonical form
than in construction.
GCD embeddings
have been computed for every $\R^3$ case with diameter $\le 100,000$.
[Indeed in 4 cases, new embeddings were found which had been overlooked ---
to considerable embarrassment on the part of its implementor --- 
by the allegedly complete search described below.]
\par
We shall need a parameterisation generating all primitive solutions
of the Diophantine equation $x^2 + y^2 + z^2 = w^2$, via the identity
\begin{align} \label{[param]}
(p^2+q^2-u^2-v^2)^2 + (2pu+2qv)^2 + (2pv-2qu)^2 = (p^2+q^2+u^2+v^2)^2. 
\end{align}
[In \cite{[Yiu01]}, where it is credited to Mordell, this same idea is employed
in an entirely different capacity --- see end of section \ref{[commalt]}.]
\par
Constraining the integer parameters $p,q,u,v$ by
\begin{eqnarray*}
  1 \le &p& \le \floor{\sqrt{w}}, \cr
  \ceil{\sqrt{\max(0, w/2-p^2)}} \le &q& \le \floor{\sqrt{w-p^2}}, \cr
  0 \le &u& \le \floor{\sqrt{\min(p^2+q^2, w-p^2-q^2)}}, \cr
  &v& = \floor{\sqrt{w-p^2-q^2-u^2}}, 
\end{eqnarray*}
before verifying $p^2 + q^2 + u^2 + v^2 = w$ will ensure that $x \le y \le z$;
though some solutions may be imprimitive or duplicated. The complete solution
requires merging scaled-up solutions for every divisor of $w$.
\par
Now the following algorithm enumerates the set of {\sl all}
essentially distinct embeddings of a given tetrahedron:
\begin{enumerate}
\item Find some rational axial pose for vertices $P,Q,R,S$ via section
\ref{[axialpose]}, noting that altitude of $S$ has ambiguous sign;
\item Generate complete solutions to $x^2 + y^2 + z^2 = u^2, v^2$ resp. 
via equation \eqref{[param]};
\item With vertex $P$ at origin, scan all possible lattice locations for
$Q,R$, retaining partial embeddings where length $QR^2 = w^2$;
\item Solve over-determined homogeneous linear equations (via null space) for
the quaternion $X$ rotating face $PQR$ from axial pose to lattice embedding;
\item If $X$ also rotates (either) $S$ to the lattice, reduce the rotated pose
$P'Q'R'S'$ to canonical form and add to result set.
\end{enumerate}
\par
A Magma implementation of this algorithm has currently found many more $\Z^3$
embeddings for tetrahedron diameter $\le 20,000$, despite evidently remaining
in need of more thorough development! A version employing assertion
\ref{[emb2emb3]} was employed (after much ill-tempered adjustment of dimension)
to verify uniqueness in $\Z^2$ for the first 1000 triangles.
\par

\section{Worked Numerical Examples}
\par
The following illustrate the complete procedure for computing a canonical
lattice embedding of a Heronian tetrahedron via quaternion GCD applied to a
rational axial pose. Throughout this section the quaternion
$s + p \i + q \j + r \k$ will be denoted $[s,p,q,r]$.
Input data comprise six edge-lengths.
\par
The first example involves the scalene case
\[ [2431,2375,1044,2296,2175,1479]. \]
Apply the currently selected permutation on vertices, say $PQRS \to QRPS$,
and via section \ref{[axialpose]} construct the corresponding rational
pose $P,Q,R,S =$ 
\[ [1,0,0,0],\ [1,1044,0,0],\ [29,18876,67925,0],\ [13,22620,8613,14616] ;\]
notice the odd scalar component (Cartesian LCD) in $R,S$.
$P,Q$ already lie in the lattice; so using section \ref{[embedR3]},
relocate $R$ by setting
\[ X = \GCD_L(R, 29) = [-5,0,0,2] \]
and transforming all four vertices by $X$ to $P,Q,R,S =$
\[ [1,0,0,0],\ [1,756,720,0],\ [1,-1144,2145,0],\ [13,10440,21837,14616] ;\]
then repeat to relocate $S$ by
\[ X = \GCD_L(S, 13) = [-2,2,2,1] \]
transforming to $P,Q,R,S =$
\[ [1,0,0,0],\ [1,396,864,432],\ [1,396,-561,2332],\ [1,1740,783,1044],\]
with every vertex now in the lattice.
Finally translate to the positive octant; scan lattice symmetries and vertex permutations as per section \ref{[canon]}, finding the canonical isomorph
\[ [1,0,0,396],\ [1,561,2332,0],\ [1,1344,1288,1740],\ [1,1425,1900,396]. \]
\par
The second example involves the same case, but permuting $PQRS \to PQSR$. 
Now the sequence follows rational pose 
\[ [1,0,0,0],\ [1,2431,0,0],\ [13,17248,24360,0],\ [17,36575,13680,10260],\]
by $X = [0,-3,-2,0]$ to
\[ [1,0,0,0],\ [1,935,2244,0],\ [1,2240,504,0],\ [17,26695,28500,-10260] ,\]
by $X = [0,3,-2,-2]$ to
\[ [1,0,0,0],\ [1,-1529,-1848,396],\ [1,-224,-1848,-1344],\ [1,-665,-2280,0] ,\]
reducing canonically to
\[ [1,0,0,396],\ [1,224,1848,1740],\ [1,665,2280,396],\ [1,1529,1848,0] ;\]
different from the previous result and so representing an essentially distinct
embedding. [In scalene cases there is no difference between strong and weak.]
For this case all 24 permutations of vertices yield one of these two results. 
\par
Nonetheless, tedious exhaustive search (section \ref{[search]}) eventually
discovers one more embedding, distinct from both above:
\[ [1,0,0,0],\ [1,224,672,2184],\ [1,665,1824,1368],\ [1,1529,1716,792], \]
showing that not every $\Z^3$ embedding is related via GCD to an axial pose.
\par
The third example employs isohedral case
\[ [8484,6625,6409,6409,6625,8484]. \]
This embedding might be thought to lack variety: whatever the vertex
permutation, any axial pose, say
\[ [1,0,0,0],\ [1,6625,0,0],\ [5,28224,31668,0],\ [5,4901,22932,21840] \]
which happens to transform in a single shot by $X = [0,-1,-2,0]$ to 
\[ [1,0,0,0],\ [1,-3975,5300,0],\ [1,1680,8316,0],\ [1,3081,3536,-4368] ,\]
always yields the same strongly canonical
\[ [1,0,0,1401],\ [1,0,3016,7056],\ [1,0,8316,3081],\ [1,4368,4780,0] .\]
Invoking only weak canonicity reduplicates this as 4 vertex-permuted variations.
More noteworthy are the results of exhaustive search in this case, yielding
in total 9(36) distinct embeddings.
\delete{
  [[1,0,0,1401], [1,0,3016,7056], [1,0,8316,3081], [1,4368,4780,0]], // GCD 
  [[[1,0,0,2088], [1,1007,3180,7812], [1,4031,4524,0], [1,7056,3024,5700]], 
  [[1,0,0,1575], [1,696,3248,7056], [1,5400,3500,0], [1,7056,1764,5943]], 
  [[1,0,0,600], [1,377,4524,5124], [1,6048,1344,6396], [1,6425,1500,0]], 
  [[1,0,0,3016], [1,2385,3180,8316], [1,3393,4524,0], [1,7728,3024,4780]], 
  [[1,0,0,2745], [1,2268,2764,8064], [1,2268,8064,4089], [1,4704,3772,0]], 
  [[1,0,0,696], [1,0,1855,7056], [1,2436,5887,0], [1,5796,4032,5400]], 
  [[1,0,0,748], [1,0,3975,6048], [1,2436,8064,1756], [1,5796,2631,0]], 
  [[1,0,0,1972], [1,1344,8064,4240], [1,1500,2745,7812], [1,4524,4089,0]]]], 

Two GCD poses with common triangle: 
  [6120,5551,1199,3825,3825,2824], 
  [7888,4732,3620,4575,3487,857], 
  [37723,34752,9125,21964,17187,16588], 
  [39360,35321,28679,32800,32800,7929]. 
GCD yields 4 poses: 
  [45100,43911,6929,34476,40544,36975], 
  [66079,63433,15204,62146,8625,11271], 
  [72336,70499,53075,46375,48841,25926], 
  [87600,79424,15184,44109,47859,39565]. 
}
\par
An intriguing feature of isohedral cases is the sparseness of their rotations:
it appears that such $X \ne 1$ have just two nonzero components,
which furthermore are very small --- eg. $X = -(\i + 2\j)$ above.
For all 11 isohedral cases with diameter $\le 20,000$, there is some vertex
permutation for which the rotations have height at most 4.
\par

\section{Heronian Pentatopes in $\R^4$} \label{[embedR4]}
\par
Consider the free pentatope $PQRST$ specified by squared edge-lengths
\begin{align*} 
\begin{tabular}{ r c c c c c c c c c c }
          edge &$QP$&$RP$&$RQ$&$SP$&$SQ$&$SR$&$TP$&$TQ$&$TR$&$TS$ \cr
 square length &  1 &  2 &  3 &  2 &  3 &  2 &  1 &  2 &  1 &  1.  
\end{tabular}
\end{align*}
$PQRST$ can be rationally posed in $\Q^4$,
with projective coordinates for its vertices:
\begin{align*}
    P &= [1,0,0,0,0], \cr
    Q &= [1,1,0,0,0], \cr
    R &= [1,0,1,1,0], \cr
    S &= [1,0,1,0,1], \cr
    T &= [3,0,1,2,2]. 
\end{align*}
Taking the determinant shows $PQRST$ is proper, with volume $1/24$.
\par
However,
\begin{Thm} \label{[counter4]}
The simplex $PQRST$ above has no lattice embedding in $\Z^4$.
\end{Thm}
\begin{proof}
We may choose lattice axes so that $T$ lies at the origin
\[ T = [1,0,0,0,0]. \]
Since $TS^2 = TR^2 = TP^2 = 1$ we may choose
\[ S = [1,0,0,0,1],\quad R = [1,0,0,1,0],\quad P = [1,1,0,0,0]. \]
Since $QP^2 = 1$, and the volume is nonzero,
the vacant column must be filled by choosing
\[ Q = [1,0,1,0,0]. \]
Now $RQ^2 = 2$, contradicting $RQ^2 = 3$ in the specification.
\end{proof}
\par
The point of all this is that assertions \ref{[poseall2]}, \ref{[poseall3]}
require only that a polytope be rationally embedded, and have edge lengths
squaring to integers, in order to be embeddable in $\Z^2$, $\Z^3$ resp.
The counter-example above demonstrates that these constraints are insufficient
to ensure embeddability in $\Z^4$; it follows that any theorem concerning
Heronian embeddability in higher dimensions must impose stronger conditions.
\par
Of course, that's assuming anything can be found to embed.
Sascha Kurz \cite{[Kur12]} recently completed the enumeration of primitive
Heronian simplices to diameter $600,000$.
He reports $41563542$ triangles, $2526$ tetrahedra,
and 0 pentatope --- that is none, nought, zero, zilch.
\par
The final paragraph of \cite{[Yiu01]} mis-states the embedding constraints.
Simplices with integer edges are discussed in \cite{[Kur10]}, \cite{[Kur11]}.
Methods for constructing Heronian simplices are discussed in detail in 
\cite{[Kur08]}, \cite{[Buc92]}, \cite{[Smi12]}.
\par

\section{Remarks on Alternative Approaches} \label{[commalt]} 
\par
During the course of an interesting but in places unspecific paper covering
many of the same topics as here, the proof of \cite{[Fri01]} theorem~4.4
develops a construction for $\Z^2$ triangle embeddings similar in spirit 
to our proof of assertion \ref{[embed2]}; however, there are significant differences.
\par
In the first place no GCD is explicitly invoked: instead the LCD $r$
is factored into rational primes $p$, then by induction each factor is
dispatched via a new rotation $X$ --- obtained in essentially the same manner
as above. From a computational viewpoint, this decomposition is both more
elaborate and slower.
\par
When it comes to dealing with the set however, instead of each point
dispatched in turn via a new rotation, every point is (partially) embedded
via the same $X$. Computationally speaking, what is lost on the swings ---
the number of factors of $q$ --- is (perhaps) regained on the roundabouts
--- the number $k$ of points.
\par
So we should enquire whether --- in $\R^3$ especially --- our induction on
$k$ might also be aggregated, winning on both rides.
Instead of transforming by $X = \GCD_L(Q, q)$ for each point $Q$ in turn,
consider just setting $X = \GCD_L(S, s)$ for the final point $S$ of the
rational pose, in the hope that $X$ will also embed every other point.
\par
Clearly this could work only if all previous $q\fac s$. But the situation
is actually rather worse: it turns out that, after rotation and reduction,
the scalar $Q'$ remains $\GCD(s/q, q)$ rather than unity. So rotated points
may miss the lattice unless, for each scalar LCD $q$ and prime factor $p$,
either $p$ divides $q$ to the same power as $p$ divides $s$, or $p$ does not
divide $q$ at all.
\par
And in practice cases do frequently occur where, besides denominators
involving various powers of the same prime, the rational pose has no LCD
equalling their LCM:
for example (for various permutations of the vertices) the early case
\[ [ 160, 153, 25, 120, 56, 39 ] \]
yields some axial poses with scalars $q = 1,1,25,65$. Tough luck!
\delete{
    [ 160, 153, 25, 120, 56, 39 ],  
    [[ 1, 0, 0, 0 ], 
    [ 1, 25, 0, 0 ], 
    [ 25, 1408, 3744, 0 ], 
    [ 65, 2912, 2016, 840 ]]. 
}
\par
Turning to $\R^3$, the paper provides a partial solution to the existence
problem: \cite{[Fri01]} theorem~4.6 asserts that computation has established
the existence of suitable rotation matrices, provided all LCDs are products
only of primes $p \le 37$. 
No details are given: one might conjecture that the method involves
some variation of the `modularisation' employed elsewhere in establishing
divisibility of volume, with consistently unreliable results.
\par
A complete contrast to the complex approach to embedding in $\R^2$ is provided
by \cite{[Yiu01]}, where the problem is reformulated in terms of a quadratic
Diophantine equation, for which the complete parametric solution \eqref{[param]}
is known. It seems improbable that an analogous attack could easily be mounted
in $\R^3$: corresponding parameterisations for simultaneous quadratic and
cubic Diophantine equations are (as far as we are aware) currently unknown,
and quite possibly do not exist.
\par
Finally, shortly after a preprint of an earlier version of this article had
been posted, we learned of \cite{[Mar12]} which coincidentally turns out
to cover much of the same ground, though via a technique closer in spirit
to that of Fricke. Swallowing our battered priority, we exchanged manuscripts
with the authors, only for them to point out an error in the proof our main
theorem; for which we must naturally thank them sincerely, albeit through
gritted teeth!
\par

\delete{

To pose a set of points axially in |Q^3 (sect. 9) it is sufficient that:
  all edge-lengths are square roots of rationals;
  one edge-length, one face area adjoining, and the volume are rational.
To then embed the set in lattice |Z^3 (sect. 10) it is sufficient that:
  all edge-lengths are square roots of integers.
In short, it suffices that all edge-lengths square to integers,
  and there is some monotonic flag having rational measures.
Whether these are also necessary I do not know!

  Define an edge to be "single" if it has odd integer length, 
    "double" if even. 
A Heronian triangle is either 
    "single" with 2 single + 1 double edge; or 
    "double with 3 double edges. 
A Heronian tetrahedron is either 
    "even" with 0 double + 4 single faces, 2 double + 4 single edges; 
    "odd" with 1 + 3 faces, 3 + 3 edges; or 
    "double" with 4 + 0 faces, 6 + 0 edges (imprimitive). 
A Heronian 4-simplex is either 
    "even", with 3 even + 2 odd solids, 4 double + 6 single edges; 
    "odd", 1 double + 4 odd solids, 6 double + 4 single edges; or 
    "double" with 5 double solids, 10 double edges (imprimitive). 
This pattern appears to continue into higher dimensions in 
an obvious fashion; though a little degenerate for dim = 1,2, 
where "single" = "odd" or "even". 

Proof of 4-space even/odd/double classification: 
Suppose pentatope has p even, q odd, r double cells: p + q + r = 5, 
Numbers of faces of each type must be even: 
  4p + 3q = 0 mod 2, q + 4r = 0 mod 2; so q = 0 or 2 or 4. 
Numbers of edges of each type must divide by 3: 
  4p + 3q = 0 mod 3, 2p + 3q + 6r = 0 mod 3; so p = 0 or 3. 
Hence p,q,r = 3,2,0; 0,4,1; 0,0,5; as claimed earlier.  
  Also solutions 3,0,2; 0,2,3; then 9 and 10 double edges, 
so would have to be double pentatope. QED 

}


\bibliographystyle{amsplain} 

\providecommand{\bysame}{\leavevmode\hbox to3em{\hrulefill}\thinspace}
\providecommand{\MR}{\relax\ifhmode\unskip\space\fi MR }
\providecommand{\MRhref}[2]{%
  \href{http://www.ams.org/mathscinet-getitem?mr=#1}{#2}
}
\providecommand{\href}[2]{#2}

\nocite{*} 

\clearpage
\end{document}